
\magnification 1200
\hsize 16.5truecm
\vsize 23truecm
\hoffset=-.4truecm
\parindent 10pt
\tolerance 1000

\font\ninerm=cmr9 scaled 1000
\font\title=cmbx10 scaled 1200
\def\trait (#1) (#2) (#3){\vrule width #1pt height #2pt depth #3pt}
\def\qed{\hfill \trait (0.1) (6) (0) \trait (6) (0.1) (0) \kern-6pt
\trait (6) (6) (-5.9) \trait (0.1) (6) (0) }

\font\sixrm=cmr6
\newcount\tagno \tagno=0                        
\newcount\thmno \thmno=0                        
\newcount\bibno \bibno=0                        
\newcount\chapno\chapno=0                       
\newcount\verno            
\newif\ifproofmode
\newif\ifwanted
\wantedfalse
\newif\ifindexed
\indexedfalse
\catcode`@=12 

\def\parag#1#2{\goodbreak\bigskip\bigskip\noindent
                   {\bf #1.\ \ #2}
                   \nobreak\bigskip}

\long\def\th#1#2{\goodbreak\bigskip\noindent
                {\bf Theorem #1.\ \ \sl #2}}
\long\def\lemma#1#2{\goodbreak\bigskip\noindent
                {\bf Lemma #1.\ \ \sl #2}}

\long\def\cor#1#2{\goodbreak\bigskip\noindent
                {\bf Corollary #1.\ \ \sl #2}}

\long\def\rem#1#2{\goodbreak\bigskip\noindent
                 {\bf Remark #1.\ \ \rm #2}}

\def\proof{\vskip.
1cm\noindent{\sl Proof.\ \ }}

\def\sqr#1#2{\vbox{
   \hrule height .#2pt
   \hbox{\vrule width .#2pt height #1pt \kern #1pt
      \vrule width .#2pt}
   \hrule height .#2pt }}
\def\square{\sqr74}

\mathchardef\emptyset="001F \mathchardef\hyphen="002D

\def\a{\alpha}
\def\b{\beta}
\def \D {{\cal D}}
\def \dps{\displaystyle}
\def \div {\mathop {\rm div}\nolimits}

\def \e{\varepsilon}

\def \H {{\cal H}}

\def \liminf {\mathop {\rm liminf} \limits}
\def \lim {\mathop {\rm lim} \limits}
\def \limsup {\mathop {\rm limsup} \limits}

\def \N {{{\rm I} \kern - .15 em {\rm N}}}
\def \Om {\Omega}
\def\bO {\overline\Om}

\def \qed {\hfill\square}
\def \res{\mathop{\hbox{\vrule height 7pt width .5pt depth 0pt
\vrule height .5pt width 6pt depth 0pt}}\nolimits}
\def \R {{{\rm I} \kern - .15 em {\rm R}}}
\def \rh {\mathop{\rightharpoonup}\limits}
\def \ren {{{\rm I} \kern - .15 em {\rm R}}^n}
\def \renn {{{\rm I} \kern - .15 em {\rm R}}^{n^2}}
\def \renns {{{\rm I} \kern - .15 em {\rm R}}^{n^2}_{{\rm sym}}}

\def \spt {\mathop {\rm spt}}
\def \square {\vrule height 1.5ex width 1.1ex depth -.4ex}
\def \sym {\mathop {\rm sym}}

\def\l{\lambda}

\def\ifundefined#1{\expandafter\ifx\csname+#1\endcsname\relax}

\def\Wanted#1{\ifundefined{#1} \wantedtrue
\immediate\write0{Wanted #1 \the\chapno.\the\thmno}\fi}

\def\Increase#1{{\global\advance#1 by 1}}
\def\Assign#1#2{\immediate
\write1{\noexpand\expandafter\noexpand\def
 \noexpand\csname+#1\endcsname{#2}}\relax
 \global\expandafter\edef\csname+#1\endcsname{#2}}

\def\pAssign#1#2{\write1{\noexpand\expandafter\noexpand\def
 \noexpand\csname+#1\endcsname{#2}}}

\def\lPut#1{\ifproofmode\llap{\hbox{\sixrm #1\ \ \ }}\fi}
\def\rPut#1{\ifproofmode$^{\hbox{\sixrm #1}}$\fi}

\def\chp#1{\global\tagno=0\global\thmno=0\Increase\chapno
\Assign{#1} {\the\chapno}{\lPut{#1}\the\chapno}}
\def\thm#1{\Increase\thmno
\Assign{#1}
{\the\chapno.\the\thmno}\the\chapno.\the\thmno\rPut{#1}}

\def\frm#1{\Increase\tagno
\Assign{#1}{\the\chapno.\the\tagno}\lPut{#1}
{\the\chapno.\the\tagno}}


\def\bib#1{\Increase\bibno
\Assign{#1}{\the\bibno}\lPut{#1}{\the\bibno}}


\def\pgp#1{\pAssign{#1/}{\the\pageno}}


\def\ix#1#2#3{\pAssign{#2}{\the\pageno}
\immediate\write#1{\noexpand\idxitem{#3}
{\noexpand\csname+#2\endcsname}}}


\def\rf#1{\Wanted{#1}\csname+#1\endcsname\relax\rPut {#1}}


\def\rfp#1{\Wanted{#1}\csname+#1/\endcsname\relax\rPut{#1}}

\Increase\verno \immediate\openout1=\jobname.aux

\immediate\write1{\noexpand\verno=\the\verno}

\ifindexed \immediate\openout2=\jobname.idx
\immediate\openout3=\jobname.sym \fi


\nopagenumbers
\headline={\toto}
\def\toto{\ifodd\pageno\rightheadline \else\leftheadline\fi}
\def\rightheadline{\ifnum \pageno=1\relax\else\ohd\fi}
\def\leftheadline{\ehd}
\voffset=2\baselineskip
\def\ohd{\centerline{\tenrm \sl \AUTHORS}}
\def\ehd{\centerline{\tenrm  \sl \TITLE}}
%
%
%
%
\def\AUTHORS{G. Bouchitt\'e}
\def\TITLE{Optimization of light structures: the
vanishing mass conjecture}


\vglue2.0truecm
\centerline{\title {Optimization of light structures: the
vanishing mass conjecture.}}


\vskip 1truecm
\centerline{\tenrm { G. B{\sevenrm OUCHITTE	} }}
\vskip3truecm
\noindent
%
{\bf Abstract:}
{\ninerm We consider the shape optimization problem
which consists in placing a given mass $m$ of  elastic material in a
  design region so that the compliance is minimal. Having in mind optimal
light structures,
our purpose is to show
that  the problem of finding the
stiffest shape configuration simplifies as the total mass $m$ tends to
zero: we  propose an explicit relaxed formulation where the compliance
appears after rescaling as a convex functional of the relative density of
mass.
 This allows us to write
necessary and sufficient optimality conditions for light structures following
the Monge-Kantorovich approach  developed recently in [5].

}
\parag{\chp{secint}} {\bf Introduction} \tagno=0 \thmno=0
 Since the beginning of the
mathematical theory of elasticity it was possible to consider from a rigorous
point of view the problem of finding the structure that, for a given system
$f$ of loads, gives the best resistance in terms of minimal compliance. In
other words, an elastic structure is optimal if the corresponding displacement
$u$ is such that the total work $\int f\cdot u\,dx$ is minimal. However, even
if the setting of the problem does not require particular mathematical tools,
only in the last two decades there has been a deep understanding of {\it shape
optimization problems} from a mathematical point of view. This was mainly due
to the dramatic improvement in the field impressed by the powerful theories of
homogenization and $\Gamma$-convergence which have been developed meanwhile.

What became clear soon was that in a large number of situations the optimal
shape does not exist, and the existence of an optimal solution must be
intended only in a {\it relaxed sense}. The form of the relaxed optimization
problem was first studied (see [22,23]) in the so
called {\it scalar case} where the physical problem only
involves state variables with value in $\R$, like the
problem of optimal mixtures of two given conductors. In this
case the relaxed solutions have been completely studied, and
identified as symmetric matrices with bounded and measurable
coefficients, whose eigenvalues satisfy some suitable {\it
bounds}. A similar result was also obtained in the
elasticity problem (see for instance [16]) for optimal
mixtures of two homogeneous and isotropic materials.
 In almost all cases which have been considered, the optimal relaxed
solution is not isotropic  and this was interpreted by saying that an
optimal shape does not exist and minimizing sequences are composed by {\it
laminates}.

We want to emphasize that the case of optimal elastic structures, or also
simply the study of optimal shapes of a given conductor, seems to have an
additional difficulty with respect to the problem of optimal mixtures. Indeed,
the first correspond to the case of optimal mixtures when one of the two
materials (or conductors) has the elasticity constants (or the conductivity
coefficient) equal to zero. In this case, due to a lack of uniform
ellipticity, it is known that among all possible relaxed problems, obtained as
limits of sequences of elliptic problems on classical domains, there are some
that are not of {\it local} type, and it is not clear if these {\it nonlocal}
relaxed solutions could be optimal. This interesting direction of research has
been developed recently in the scalar case in [2] and [20]
showing deep
connections with the theory of Dirichlet forms, and more recently in the
case of elasticity in [14].

\bigskip
Here we are interested in a apparently more difficult problem which
consists in finding
the asymptotic of the previous shape optimization problem when the total
volume tends to zero.
In other words we are trying to give a mathematical
fundation  for what we call optimal light structures.
In [1] it has been proposed to solve this problem by
 following three steps: 1) describe optimal mixtures of
two elastic materials (homogenization) 2) Pass to the limit when the
rigidity constants of the weak material
tend to zero, 3) Pass to the limit as the ratio of void goes to $1$ (high
porosity limit).

However I got the impression that this approach is very
heavy and, as far as I know, in view of the difficulties explained above no
real
mathematical justification has been given until now for what concerns steps
2 and 3. In addition
the occurrence of concentration on lower dimensional structures  expected in
many cases by engineers and
manufacturers pushes a priori for searching optimality out of the class of
microstructures.

The aim of this paper is to propose a new direct approach leading to a very
simple formula
for the limit compliance where the light structure is described in term of
the  density distribution of
material. This density is a possibly concentrated non negative measure and
the corresponding energy functional
turns out to be {\it convex}.  As a consequence we may use the framework I
recently developed
in collaboration with G.Buttazzo and P. Seppecher [4,5,6] which allows us to see the
optimal measure as the multiplier of a linear programming problem and also
to write necessary and
sufficient optimality conditions  (Monge-Kantorovich system).

The plan of the paper is the following. In Section 2 we present the
rescaled shape
optimization problem associated with a small total mass $\e$. It is written
as the minimization  of a
functional
$c_\e(\mu)$ defined on probability measures. Then
we state the form of the
$\Gamma$-limit of the sequence  $\{c_\e\}$  as $\e\to 0$. The model of
Michell truss like structures is
recovered as a particular case in dimension 2. In section 3, we introduce a
framework suitable for dealing
with lower dimensional structures in $\R^n$. We derive a relaxed compliance
for structures whose dimension
is prescribed to be less than or equal to $k<n$. Applying this result for
$k=n-1$ allows us to prove the
upperbound inequality for the $\Gamma$-limit of $\{c_\e\}$. In the last
section,
the lower bound inequality is
presented as a consequence of what we call the  {\it vanishing mass
conjecture}. Some geometrical
arguments are given to support this conjecture.

\parag{\chp{secint}} {\bf Setting of the problem and the vanishing mass
model.} \tagno=0 \thmno=0
\noindent{\bf Notations.}\ In what follows $\Omega$ is a bounded Lipschitz
connected open subset of $\R^n$ ,
$\Sigma$ is a compact subset of $\bO$ and $F$ is an element of
${\cal M}(\bO;\R^n)$, the class of all $\R^n$-valued measures on $\bO$ with
finite total variation.

The class of smooth
displacements we consider is the Schwartz space $\D:=\D(\R^n;\R^n)$ of
$C^\infty$
functions with compact support; similarly, the notation $\D'(\R^n;\R^n)$ stands
for the space of vector valued distributions and, for a given nonnegative
measure
$\mu$, $L^2_\mu(\R^n;\R^d)$ denotes the space of $p$-integrable functions.
 The symbol $\cdot$ stands
for the Euclidean scalar product between two vectors in $\ren$ or
between two matrices in $\renn$.
\bigskip
The elastic structure is placed in $\Om$ and occupies a region (an open subset)
$\omega\subset\Om$ of prescribed volume $m$. It is clamped
on the part of $\omega$  in contact with $\Sigma$ and
 it has to support the given load $F$. The (possibly infinite) compliance
associated with this configuration is given by
$$  c(\omega,F,\Sigma):= - \inf \left\{ \int_\omega j(e(u)) dx - <F,u> \ :\
u\in  \D, \ u=0 \
\hbox{on
$\Sigma$}\right\}\ ,\leqno(\frm{comp})$$
where $j: \renns \mapsto [0,+\infty)$ is a quadratic form characterizing
the elastic properties of the material
and $e(u)$ denotes the symmetrized tensor of deformations i.e. $e(u)_{i j}=
{1\over 2}({\partial u_i\over
\partial x_j }+ {\partial u_j\over \partial x_i})$.
In the 3D case a model example is given by the isotropic  elasticity with
Lam\'e coefficients $\a,\b$
($3\a +2\b >0$):
$$ j(z) = {1\over 2} \a  |tr z|^2  + \b |z|^2\ .\leqno(\frm{lame})$$
A classical situation is when $\Sigma $ is a part of the boundary of $\Om$
and $F$ is of the kind
$$ <F,u> = \int_\Om f u dx + \int_{\partial\Om\setminus \Sigma}  g u \, d
{\cal H}^{n-1}\qquad \hbox{being}
\ f\in L^2(\Om;\R^3)\ ,\ g\in L^2(\partial\Om\setminus \Sigma;\R^3).
$$
In this case and if the data $f,g$ are compatible (their supports need to be
contained  in ${\overline
\omega}$), the  infimum in (\rf{comp}) is finite and minimizers can be
searched in the
Sobolev space
$W^{1,2}(\omega;\R^n)$.
Let us stress the fact that much more general situations will be considered
in our framework since we
intend further to consider concentrated loads (for example Dirac delta) as well
as sets $\omega$ whose measure becomes very small.

Our interest is to pass to the limit as $m\to 0$ in the following
variational problem
$$  \inf \{ c(\omega) \ ;\  |\omega|=m\}\ .\leqno(\frm{optim})$$
As the data $\Sigma, F$ are kept fixed all along the paper we will write
$c(\omega,F,\Sigma)=c(\omega)$.

\bigskip
\noindent{\bf Rescaled Problem.\ } We set the total mass $m$ to be a small
parameter $\e$ tending to $0$. It is
easy to check that the infimum in (\rf{optim}) scales as ${1\over\e}$ as
$\e\to 0$. Indeed  by making the change of
variables
$v=\e u$ in the integral, we obtain the identity:
$$ \e \, c(\omega)\, = \, c_\e (\mu) \quad,\qquad  \mu(x)={1\over \e} \,
1_\omega(x)\ ,$$
where $1_\omega$ denotes the characteristic function of $\omega$ and $c_\e$
is the functional on ${\cal
M}_+(\bO)$ defined by:
$$ c_\e(\mu) :=\left\{\eqalign{
& - \inf\left\{\int j(e(u)  d\mu - <F,u> :\ u\in \D, u=0 \hbox{\ on
$\Sigma$} \right\}\cr
 &\hbox{ if $\displaystyle\mu= \mu(x) \, dx\  ,  \mu(x)\in
\big\{0,{1\over\e}\big\}$ }
\quad  (+\infty  \ \hbox{ otherwise})\ .
}\right.
 \leqno(\frm{ce})$$
The total mass constraint becomes $\int_{\bO} d\mu =1$ so that (\rf{optim})
can be restated for $m=\e$ and
after rescaling as the minimization of $c_\e(\mu)$ over all probability
measures on $\mu$ on $\bO$.
It is clear that the latter infimum decreases if we relax the constraint on
the density of $\mu$ which appears
in (\rf{ce}). More precisely, we have
$$  \inf \{c_\e(\mu) \ :\  \int_\Om d\mu=1 \}  \ \ge\  \inf \{c(\mu) \ :\
\int_{\bO} d\mu=1 \}\
,\leqno(\frm{gap})$$ where $c(\mu)$ is defined for every $\mu\in{\cal
M}_+(\bO)$ by
$$ c(\mu):=- \inf\left\{\int j(e(u)  d\mu - <F,u> :\ u\in \D, u=0 \hbox{\
on $\Sigma$} \right\}\
.\leqno(\frm{mass})$$
We notice that the functional $c(\mu)$ defined above has the form already
used in  [5] to modelize
mass optimization problems. It is a convex lower semicontinous functional
on measures. Unfortunately, as we
will see later,   a gap in (\rf{gap}) will occur in general except
particular situations where the
original mechanical problem can be reduced to a scalar setting.

Now the natural procedure in order to pass to the limit as $\e\to 0$ in the
left hand side of (\rf{gap})
consists in computing the $\Gamma$-limit of $c_\e$ with respect to the
weak (star) convergence of measures
on the compact $\bO$. Let us introduce the following integrand:
$$ \bar j(z) := \sup \left\{ z \cdot \xi - j^*(\xi) \ : \ \xi\in \renns  \
,\ det\, \xi=0 \right\}\ ,
\leqno(\frm{new})$$
where $j^*$ denotes the Fenchel transform of $\xi^*$ (i.e.  $j^*(\xi)=
\sup\{ z\cdot \xi - j(z) \ :\
z\in\renns\}$). By Lemma \rf{jkk} in section 3, it turns out that $\bar j$
is  convex continuous and satisfies
$$   \bar j (z)  \le j(z)  \quad  \hbox{for all $z$} \quad,\quad
(\bar j)^*(\xi) = j^*(\xi) \quad \hbox{whenever $rank (\xi)<n$ }\ .$$

Our claim is two fold: first we say that the $\Gamma$ limit of the sequence
$\{c_\e\}$ exists and is a
convex functional; then we  claim that this limit can be represented like
in (\rf{mass}) but substituting $j$
with the new integrand $\bar j$.
Precisely we conjecture the following result:

\th{\thm{conj1}} ({\sl Conjecture }) {\it  The $\Gamma$-limit of $c_\e$
with respect to the weak convergence of
measures on $\bO$ is the functional $E(\mu)$ defined on ${\cal M}_+(\bO)$ by
$$ E(\mu):=- \inf\left\{\int {\bar j}(e(u)  d\mu - <F,u> :\ u\in \D, u=0
\hbox{\ on $\Sigma$} \right\}\
.\leqno(\frm{formula})$$  }
To prove this theorem, we need to show that:

a) For every sequence $\{\mu_\e\}$ such that $\mu_\e \rh \mu$, there holds
$\liminf_\e c_\e(\mu_\e) \ge E(\mu)$.

b)  For every measure $\mu \in {\cal M}_+(\bO)$, there exists a sequence
$\{\mu_\e\}$ such that
$\mu_\e \rh\mu$ and $\limsup_\e c_\e(\mu_\e) \le E(\mu)$.

The proof of b) will be sketched in section 3 where structures of lower
dimension are considered.
The proof of a) will be straightforward in the case where $\mu$ is 
concentrated
on lower dimensional manifolds.
The general case will be seen as a consequence of geometrical properties
related with sets of vanishing measure
( in section 4, we call it  {\it vanishing mass conjecture}).

Let us  now compute $\bar j$ for particular  functions $j$ given in the form
(\rf{lame}).

a) Assume $\a=0, \b={1\over 2}$ and $n=3$. Then, we can exploit formula
(\rf{new}) writing the symmetric tensor
$\xi$ in the form
$ \xi =  \tau_1 e\otimes e  + \tau_2 e^\perp\otimes e^\perp$, where $e$ is
a unit vector and  $\tau_1,\tau_2,0$
are  the eigenvalues of $\xi$. We obtain
$$\eqalign{
{\bar j}(z) &= \sup\left\{ \tau_1 (ze\cdot e) +\tau_2 (ze^\perp\cdot
e^\perp)  - {1\over 2} (\tau_1^2 +\tau_2^2)
\ :\ |e|=1,\ \tau_1,\tau_2\in\R\right\}   \cr
 &= \sup\left\{ {1\over 2} \left((ze\cdot e)^2 +(ze^\perp\cdot
e^\perp)^2\right) \ :\ |e|=1 \right\} \cr
&= {1\over 2} ( \l_1(z)^2 + \l_2(z)^2)\ ,
\cr }$$
where $|\l_1(z)| \ge |\l_2(z)| \ge |\l_3(z)|$ are the eigenvalues of the
tensor $z$.
\medskip
The computation of $({\bar j})^*$ is rather complicated and of course can be
generalized for any pair of Lam\'e
coefficients $\a,\b$. In fact  we recover this way the formulae obtained in
[1] where explicit form
of the relaxed stress potential are given (these formulae turn out to be in
agreement with  our $({\bar j})^*$).

b)  The case $n=2$ is simpler. Setting $\gamma= {\a+2\b \over 4 \b(\a+\b)}
$, we have that
$ j^*(\xi)= {1\over 2} \gamma \tau^2$ holds for any rank one tensor
of the kind $\xi=\tau e\otimes e$ where  $|e|=1$.
Then denoting by
$|\l_1(z)|
\ge |\l_2(z)|$ the eigenvalues of
$z\in
\R^4_{\sym}$ and by $\tau_1(\xi),\tau_2(\xi)$ the eigenvalues of $\xi$
 we find easily
$$ {\bar j}(z) = {1\over 2\gamma}\  |\l_1(z)|^2\quad,\quad ({\bar j})^*(\xi)=
{1\over 2} \gamma\
(|\tau_1(\xi)| + |\tau_2(\xi)|)^2  .\leqno(\frm{mich})$$

We notice that in both examples the new potential $\bar j$ is not quadratic any
more. However it remains always convex and homogeneous of degree 2. Accordingly we
are able
to treat the minimization of the compliance
relative to $\bar j$ using many tools developed in [5].
This is summarized in the following corollary. It is convenient to
introduce the following convex continuous
positively 1-homogeneous integrands on $\renns$:
$$ \rho(z):= \inf \{ t>0 :{\bar  j}({z\over t}) \le {1\over 2} \}\quad,\quad
\rho^0 (\xi):= \sup \{ \xi \cdot z \ :\
{\bar j}(z)
\le {1\over 2}
\}\ .$$
As $\bar j$ is $2$-homogeneous, we have ${\bar j}(z)={1\over 2}\,
\rho(z)^2$ and
${\bar j}^*(\xi)={1\over 2}\, (\rho^0(\xi))^2$
An important quantity associated with the data $\Omega, F,\Sigma$ is given by:
$$ I(F,\Om, \Sigma):= \sup \{ <F,u>\ :\ \rho(e(u))\le 1\ , \ u\in
U^\infty(\Om;\ren),\ \  u=0 \ \hbox{on}
\Sigma\}\ ,\leqno(\frm{cons})$$
where  $ U^\infty(\Om;\ren)$ denotes the space of functions $u\in
L^\infty(\Om;\ren)$ such that $e(u)\in L^\infty(\ren;\renns)$. In   [5] (see also
[6]) it is proved that the latter supremum is
achieved. Notice that in general functions in $U^\infty$ are not Lipschitz
due to the lack of Korn's inequality
in $W^{1,\infty}$.

\cor{\thm{mk}}\ {\it
Let $\{\omega_\e\}$ a sequence of domains such that
$$c(\omega_\e) \le \inf\{ c(\omega): | \omega|=\e\} + O(\e)\qquad
\hbox{(c.f. (\rf{optim}))}\ .$$
 Then the  sequence $\mu_\e={1_{\omega_\e} dx \over \e}$ converges weakly
(up to a subsequence) to a probability
measure $\mu$ on $\bO$.  The following assertions hold:
\medskip\noindent i)\ $\mu$ solves the minimum problem:
 $$\inf \left\{ E(\mu)\ :\  \spt \mu\subset \bO,\ \int d\mu=1\right\}
\qquad \hbox{( $F$ given by
(\rf{formula}))}\ .\leqno(\frm{mass})$$
\medskip\noindent ii)\  We have:
$$   \min  \hbox{(\rf{mass})} = { \big( I(F,\Om,\Sigma)\big)^2\over 2}\ .$$
\medskip\noindent ii)\ The following equality holds:
$$ I(F,\Om,\Sigma)  = \min \left\{ \int \rho^0(\lambda) \ :\ \l\in {\cal
M}(\ren;\renns)\ ,\ - \div \l = F \ \hbox{on $\ren\setminus \Sigma$}
\right\}\ .
\leqno(\frm{dual})$$
Moreover if a vector measure $\l$ is optimal for (\rf{dual}), then $\mu=
I(F\!,\!\Om\!,\!\Sigma)\  \rho^0(\l)$ is a minimizer
of problem (\rf{mass}).
}

\rem{\thm{michell1}}
The measure $\l$ appearing in (\rf{dual}) is the stress field associated
with the equilibrium of
the optimal structure. This stress can be written in the form $\l = \sigma
\mu$ where the density $\sigma$
satisfies
$\rho^0(\sigma)$ is constant along the optimal structure represented by
the measure $\mu$. This fact is the
mathematical counterpart of some fact which is well known to engineers.
Recall here that the notation $\rho^0(\l)$ denotes the non negative
measure of density $\rho^0(\sigma)$ with
repect to $\mu$. By the 1-homogeneity of $\rho^0$  this measure is
independent of the decomposition
$\l=\sigma\mu$ (see [17]).
 In the case $n=2$, owing to (\rf{mich}) we obtain
$\rho^0(\xi)= {\sqrt \gamma}  (|\tau_1(\xi)| +|\tau_2(\xi)|)$  and the
infimum problem (\rf{dual})
becomes nothing else but the celebrated Michell's problem [17].

\rem{\thm{michell2}}  The optimal solution $\mu$ of (\rf{mass}) can be
characterized by a system of
optimality conditions (Monge-Kantorovich system, see [5] )
which involves a notion of $\mu$-tangential
derivative. It turns out that in general the solution $\mu$ is not unique and
 is very sensitive to the form of the integrand ${\bar j}$ or
equivalently   to the form of the convex set of
matrices
$\{z\in \renns: \rho(z)\le 1\}$.
In fact using the optimality conditions, it has been proved in [5,
example 5.1] that for some 2D
configurations, the minimum (\rf{mass}) can be strictly greater than the one
obtained keeping the initial
elastic potential $j(z)= |z|^2$ instead of ${\bar j}(z)=\l_1(z)^2$.
Moreover the topology of the solution
changes drastically: the numerical solution for $j$ found in [18]
gives a two dimensional positive
density (no holes) whereas for $\bar j$ (Michell's problem) a lot of
solutions made with junction of bars
can be found.

\rem{\thm{michell3}}  The same problem can be handled from the beginning in
the scalar case (heat
equation), meaning that in (\rf{comp}), $u:\ren \to \R$ , $j:\ren \to \R$
and $F\in{\cal M}(\bO;\R)$.
In this case, the situation becomes much simpler: Theorem \rf{conj1} holds
with $E(\mu)\equiv c(\mu)$. In
other words
${\bar j}=j$ and the gap appearing in (\rf{gap}) goes to zero as $\e\to 0$.
It turns out that in this case as $ \e\to 0$, minimizing structures either
concentrate on lower
dimensional manifolds or spread in many thin layers which are parallel to
the direction of the gradient of
some optimal $u$ related to problem (\rf{cons}). Let us finally notice that
in the case $\Sigma=\emptyset$
and $\int F=0$, the supremum in (\rf{cons}) is nothing else but the Monge
-Kantorovich norm distance
between the positive and negative parts of $F$.

\parag{\chp{secint}} {\bf Compliance of lower dimensional structures.}
\tagno=0 \thmno=0

These lower dimensional structure can be justified from two sides:

- they can be seen as limits of n-dimensional structures supported on sets
of vanishing measure
(this is the essence of the classical fattening approach);

- the designer can be interested in finding stiff structures made
exclusively with beams (1D structures) or
with plates (2D structures) or also with a blend of 1D-2D structures
(excluding volumic parts).
In other words, the competitors $\mu$ in the minimization problem with
respect to $c(\mu)$ given in
(\rf{mass}) have to be searched in the subclass of measures supported by
subsets of dimension
$k$ with $k=1$ or $k=2$ or $k\le 2$. Of course this additional constraint
will increase the value of the
infimum.

We introduce for every value of the integer $k$ ($k\in[1,n]$) the following
functional on ${\cal M}(\bO)$:
$$ c_k(\mu)\   :=\cases{ c(\mu) & if $\dim T_\mu(x) = k \
\mu$-a.e.
\cr +\infty & otherwise}\
\leqno(\frm{ck})$$
where  $\dim T_\mu(x)$ represents the dimensional of the tangent space to
$\mu$ at $x$.
This notion will be made precise later.

We introduce also the following convex integrand on $\renns$:
$$ j_k(z) := \sup \left\{ z \cdot z' - j^*(\xi) \ : \ \xi\in \renns  \ ,\
{\rm rank}\, \xi\le k  \right\}\
.
\leqno(\frm{jk})
$$
The properties of the $j_k$'s are summarized in the following lemma
\lemma{\thm{jkk}}\ {\it

i)  For every $k$, $j_k$ is convex, homogeneous of degree two and we have
$$0 = j_0 \le j_1 \le \dots \le
j_{n-1}\le j_n=j\quad \hbox{and}\quad j_{n-1}=\bar j\ .$$

ii) The Fenchel conjugate of $j_k$ is given by
$$ j_k^*(\xi) = \inf \left\{\int j^* d\nu \ :
\spt \nu \subset G_k,\ [\nu]=\xi
\right\}
\ ,
\leqno(\frm{fenchel})$$
where the infimum is taken on probability measures on $\renns$, $[\nu]$ is
the barycenter of $\nu$ and
$G_k$ denotes the tensors of rank not greater that $k$.

iii) We have $j_k^*(\xi) = j^*(\xi)$ whenever $\xi\in G_k$.
}
\proof i) is trivial and iii) is a consequence of ii) noticing that
$\nu=\delta_{\xi}$ is admissible when
$\xi$ belongs to $G_k$. Let us denote now by $g_k(\xi)$ the infimum which
appears in the right hand side
of (\rf{fenchel}). Clearly $g_k$ is convex and since any symmetric tensor
$\xi\in\renns$ can be decomposed
in a convex combination of at most $n$ rank one tensors, it is finite and
continuous.
Therefore proving the equality $j_k^*= g_k$ is equivalent to showing that
$g_k^*=j_k$.
For every $z\in \renns$, we have
$$\eqalign{ g_k^*(z) &= \sup_{\xi,\nu} \left\{ z\cdot\xi - \int j^* d\nu\
: \ [\nu|=\xi,\ \spt \nu\subset
G_k
\right\}
\cr
&=  \sup_{\nu} \left\{ \int \big(z\cdot s - j^*(s)\big) d\nu\ : \spt
\nu\subset G_k \right\}\
\le  j_k(z) \ .
}\ ,$$
where the last inequality becomes an equality if we choose competitors
$\nu$ of the kind $\nu=\delta_\xi$
where $\xi$ runs over $G_k$.
\qed
\bigskip

The main result of this section is the following:
\th{\thm{relaxk}} {\it Let $\overline c_k$ denote the lower semicontinous
envelope of $c_k$. Then there holds for every measure $\mu\in {\cal M}_+(\bO)$:
$$  \overline c_k(\mu)  =- \inf\left\{\int j_k(e(u)  d\mu - <F,u> :\ u\in
\D, u=0
\hbox{\ on $\Sigma$} \right\}\ .\leqno(\frm{formula})$$  }

\rem{\thm{eifel}}\ We stress that the domain of $ \overline {c_k} $ extends
to all the space ${\cal
M}_+(\bO)$. This means that the dimensional constraint is not closed;
however the approximation of higher
dimension structures (say of dimension $l$) involves some additional energy
which
corresponds to the gap between the
integrands $j_k$ and $j_l$ (see lemma \rf{jkk}). Roughly speaking a structure
obtained by using micro-structures made of beams might look like a  3D
elastic structure. However its
compliance cannot be predicted by using the original model of 3D elasticity.

\rem{\thm{eifel2}}\ The relaxation procedure for $c_k$ does not change if
we replace the equality
constraint $\dim T_\mu(x) = k \
\mu$-a.e. by the inequality constraint constraint $\dim T_\mu(x) \le k $.
The lower semicontinuous envelope of this new functional $c_k$ will be
still given by the right hand side
of (\rf{formula}). On the other hand, the identity  $\bar j=j_{n-1}$ (see
assertion i) of Lemma \rf{jkk})
implies that we have, for every $\mu$:
$$  E(\mu) \ =\ \overline {c_{n-1} }(\mu)\ .\leqno(\frm{identite})$$

\cor{\thm{limsup}}\ {\it For every probability measure $\mu$ on $\bO$,
there exists a sequence of subsets
$A_\e\subset\Om$ such that $\mu_\e \rh \mu$ and $\limsup c_\e(\mu_\e) \le
E(\mu) \ .$}

\proof (sketch) Take $\mu$ of the kind $\mu= \theta \H^{n-1}\res S$ where
$S$ is a smooth
$n-1$-dimensional manifold and $\theta$ is a continuous and positive weight
on $S$. Then using a standard
fattening method, it is possible to construct a sequence of sets $A_\e$ such
that $|A_\e|=\e$,
$\mu_\e:={1_{A_\e}\over\e} \rh \mu$ and such that $\limsup_\e c_\e(\mu_\e)
=\limsup_\e c(\mu_\e) \le
c(\mu)$. Therefore the $\Gamma$-limsup of $c_\e$ (denote it $E_+$) which is
weakly lower semicontinous
satisfies for every such
$n-1$ dimensional measure $\mu$ the inequality
$$ E_+( \mu) \le  c(\mu) = c_{n-1} (\mu)\ .$$
The previous inequality is then extended to all $n-1$ dimensional measures
yielding the inequality
$E_+ \le c_{n-1}$. The conclusion follows by passing to the lower
semicontinuous envelopes taking into account
(\rf{identite})\ . \qed

\medskip
The end of this section is devoted to a sketch of the proof of  Theorem
\rf{relaxk}. (a complete version
will be found soon in [3]). Before we need
 to give a precise meaning to the notion of $k$ dimensional structure and
for that  we make use of the
concept of tangent space $T_\mu(x)$ to a measure introduced in
[10] which  makes sense for {\it any} positive Borel measure $\mu$ on
$\ren$.  The underlying idea is to identify
every subset $S$
of $\ren$ having Hausdorff dimension $k $  with  the overlying
measure $\H ^k$, possibly weighted by a positive density
$\theta$; more in general, a multijunction made by the union of
sets $S_i$ with different dimensions $k _i$ may be described
through a positive measure $\mu$ of the kind $\sum _i \theta _i
\H^{k_i} \res S_i$.
We refer to the  papers [4, 7, 8, 9, 10, 11, 12, 13]
where this framework has been developed with  many applications in
elasticity, shape optimization and
homogenization.

\noindent
Here we only sketch the features which are useful for the understanding of
some arguments developed later.

\noindent{\bf Tangent space }. It is a $\mu$- measurable multifunction
$T_\mu(x)$ from $\ren$ into the
linear subspaces of $\ren$. The shortest way to define it (perhaps not
the more intuitive one) is the following:
consider the operator $B: (L^2_\mu)^n \mapsto L^2_\mu$ defined by
$$\cases{D(B):= \Big \{ \sigma \in (L^2 _\mu) ^n\ :\ \exists
\, C>0 \hbox{ such that } \Big |\dps{\int} \sigma \cdot \nabla u
\, d \mu \Big |\leq C \|u \|_{2, \mu} \ \forall u \in {\cal D}
\Big \} & \cr
\noalign{\medskip} B  \sigma = v\ \Longleftrightarrow \ -\div(\sigma \mu) = v
\mu \ .& \cr}
$$
Then it can be proved that the closure of $D(B)$ coincide with the set of
selections of a
unique (up to the $\mu$ a.e. equivalence class) of a multifunction. This
multifunction denoted $T_\mu(x)$
(our tangent space) is characterized by the following equality
$$  {\overline  D(B)} = \left\{ \sigma\in L^2_\mu(\ren;\ren)\ :\
\sigma(x)\in T_\mu(x)\  \mu
\hbox{a.e.}\right\}\ .\leqno(\frm{B1})$$ In what follows $P_\mu(x)$ (possibly
identified to an element of $\renns$) will
denote the orthogonal projector on
$T_\mu(x)$.

\medskip\noindent \noindent{\bf Tangential gradient }. A more intuitive
path to reach the definition of
$T_\mu(x)$ is motivated by the following problem: given a sequence
$\{u_h\}\subset \D$ such that
$(u_h,\nabla u_h) \to (u,\chi)$ in $L^2_\mu$, what can we say of the
relation between $u$ and $\chi$ ?
In other words,  we are looking for a characterization of the closure of
the set
$$  G := \{ u,\nabla u) \ :\ u\in \D(\ren) \} \ .$$
It turns out that we have
$$ (u,\chi)\in {\overline G}  \quad \Longleftrightarrow \ u\in H^1_\mu
\quad \hbox{and}\quad\exists \xi\in
L^2_\mu(\ren;T_\mu^\perp)\ :  \chi= \nabla_\mu u + \xi\ .
 \leqno(\frm{struct})$$
Here $\nabla_\mu u$ ({\it $\mu$-tangential gradient}) is defined for smooth
$u$ by
setting $\nabla_\mu u(x)=
P_\mu (\nabla u(x))$ and is extended in a unique way (as a closable
operator)  to all functions of the
Sobolev space $H^1_\mu$. This  space $H^1_\mu$ is the completion
of smooth function with respect
to the Hilbert norm $\Vert u\Vert := (\int (u^2 + |\nabla_\mu u|^2)
d\mu)^{1\over 2}$

\medskip\noindent \noindent{\bf Tangential strain }. The same scheme can be
used substituting $\nabla u$
by the strain $e(u)$ of a vector function $u\in \D(\ren;\ren)$ (see
[5, 13]. We now consider:
$$G:=\left\{ (u, e(u)) \in (L^2_\mu)^n \times L^2_\mu (\ren;\renns) \ :\
u\in \D(\ren;\ren)\right\}\ .$$
The
analogous of statement (\rf{struct}) reads as
$$ (u,\chi)\in {\overline G}  \quad \Longleftrightarrow \ u\in \D^{1,2}_\mu
\quad \hbox{and}\quad\exists
\xi\in L^2_\mu(\ren;M_\mu^\perp)\ :  \chi= e_\mu(u) + \xi\ .
 \leqno(\frm{struct})$$
Here $M_\mu(x)$ is a multifunction from $\ren$ to vector subspaces of
$\renns$ which can be defined using
the analogous of operator $B$ defined in (\rf{B1}) where now $\sigma\in
L^2_\mu(\ren;\renns)$.
In [12], it is shown that under very mild regularity assumptions on $\mu$
(regularity by blow-up), the
relation between $M_\mu$ and $T_\mu$ is explicit:
$$  M_\mu(x) = \left\{ P_\mu(x) \xi P_\mu(x) \ :\ \xi\in \renns \right\}\ .
\leqno(\frm{Mmu})$$
Denoting by $Q_\mu(x)$ the orthogonal projector on $M_\mu(x)$, the
tangential strain is  defined for
$C^1$ functions by $e_\mu(x)= [Q_\mu(x)](e(u)(x))$ extended by continuity
to $\D^{1,2}_\mu$ the completion  of
the space of smooth deformations with respect to the $L^2_\mu$-energy.

\medskip\noindent \noindent{\bf Stress formulation for the compliance }\ It
is a matter of classical convex
analysis to show that $c(\mu)$ given by (\rf{mass}) can also be written as
$$ c(\mu) = \min \left\{ \int j^*(\sigma) \, d\mu\ :\ \sigma\in
L^2_\mu(\ren;\renns)\ ,\ -\div (\sigma
\mu)= F \quad \hbox{on}\ \ren\setminus \Sigma
\right\}\ .
\leqno(\frm{stress})$$
We stress the fact that the competitors $\sigma$ in (\rf{stress}) belong to
$D(B)$ and therefore satisfies
$\sigma(x)\in M_\mu(x)$, $\mu\, $ a .e. In particular if $\mu$ is
associated with a $k$-dimensional
structure, then we have ${\rm rank}\ \sigma(x) \le k$  $\mu\, $ a .e

We notice also that a similar representation formula for
$\overline{c_k}(\mu)$,$E(\mu)$ are obtained by
simply substituting $j^*$ respectively with $j_k^*$ and $({\bar j})^*$ in
(\rf{stress}).

\bigskip
\noindent {\bf Proof of
Theorem \rf{relaxk}  (sketch):}

i) {\sl Lowerbound:} \ Denote by $E_k(\mu)$ the right hand side of
(\rf{formula}).
It is a convex l.s.c. functional and we have
$$ E_k(\mu)  = \min \left\{ \int j_k^*(\sigma) \, d\mu\ :\ \sigma\in
L^2_\mu(\ren;\renns)\ ,\ -\div (\sigma
\mu)= F \quad \hbox{on}\ \ren\setminus \Sigma
\right\}\ .$$
Assume  that $\dim (T_\mu(x)) \le k$. Then arguing as above, the competitors
$\sigma$ have rank $\le k$ and
therefore $j_k^*(\sigma)=j^*(\sigma)$. Owing to (\rf{stress}), we deduce
that, for such measures $\mu$, there holds
$E_k(\mu)= c(\mu)=c_k(\mu)$. Thus the inequality  $c_k \ge E_k$ holds over all
${\cal M}_+(\bO)$.
Passing to the lower semicontinuous envelopes, we
deduce that $\overline{c_k} \ge E_k$.

ii) {\sl Upperbound:}\ We restrict here to the case where $\mu$ is the
Lebesgue measure. For any $Y$-periodic
measure $\nu$ of dimension $\le k$ (here $Y$ is the possibly rotated unit
cube in $\ren$ and we normalize $\nu$ so that $\nu(Y)=1$), we consider
$\mu_\e:= \nu({x\over\e})$ which clearly converges to $\mu$. The limit of
$c(\mu_\e)$ can be predicted from
the theory of homogenization for thin structures which has been developed
in several directions (see [15] and
also [8, 11, 12]). The upperbound inequality for
$\overline{c_k}$ can be deduced by
optimizing other such periodic measures $\nu$ using in particular the
homogenized stress potential
representation  obtained in [12] together with an argument of
localization which allows to commute infimum
and integral.

\qed

\parag{\chp{secint}} {\bf Lowerbound inequality and the vanishing mass
conjecture.} \tagno=0 \thmno=0
To complete the proof of Theorem \rf{conj1}, we need to show that for every
sequence of sets $\{A_\e\}$
 such that $|A_\e|=\e$, the following implication holds:
$$ \mu_\e :={1_{A_\e}\over \e} \rh \mu\quad \Rightarrow\quad  \liminf_\e
c_\e(\mu_\e) \ge E(\mu)\
.\leqno(\frm{liminf})$$
Note that in general $c(\mu) \le E(\mu)$ (see (\rf{formula}) whereas
the lower semicontinuity of $c(\mu)$ implies only the inequality
$\liminf_\e c_\e(\mu_\e) =\liminf_\e c(\mu_\e) \ge c(\mu) \ .$

In fact the inequality (\rf{liminf}) is straightforward in the case where
$\mu$ has no $n$-dimensional
part. More precisely, we have
\lemma{\thm{low}} {\it  Let $\mu\in {\cal M}_+(\bO)$ such that $\dim
T_\mu(x)<n\ ,\ \mu\, $a.e.
Then
$$\Gamma-\lim_\e c_\e (\mu) = E(\mu) \ .$$
}
\proof\ According to the observation made above, we are done if we show
that such measures satisfy
$E(\mu)=c(\mu)$. We argue on the stress formulations of $c(\mu)$ and
$E(\mu)$ (see (\rf{stress}), simply
noticing that all admissible stress fields
$\sigma$  in the minimum problem (\rf{stress})
 satisfy $\sigma\in M_\mu(x)$   and therefore ${\rm rank}\, \sigma(x)\le
\dim (T_\mu(x))\le n-1\ \mu\,
$a.e. Thus by Lemma \rf{jkk}, we have $({\bar j})^*(\sigma)=
j_{n-1}^*(\sigma)=j^*(\sigma)$ yielding
the equality $E(\mu)=c(\mu)$.

\qed

The approach we suggest for attacking the general case consists, like in
previous
lemma, in considering the dual
formulations. Let $\sigma_\e$ be the optimal stress  related to $\mu_\e\,
$; it satisfies
$$ c_\e(\mu_\e) = \int j^*(\sigma_\e) d\mu_\e \quad ,\quad -\div \sigma_\e
\, \mu_\e = F \quad \hbox{on}\
\ren\setminus \Sigma\ .$$
Assuming that $c_\e(\mu_\e)$ is uniformly bounded (otherwise (\rf{liminf})
is trivial), by the growth
condition on $j$ we have that
 $$\sup_\e \int |\sigma_\e|^2 \, d\mu_\e \ <\ +\infty \ .\leqno(\frm{bound})$$
From (\rf{bound}) it follows that possibly passing to subsequences, there
exists a suitable $\sigma\in
L^2_\mu(\ren;\renns)$ such that $\sigma_\e \,\mu_\e \rh \sigma\, \mu$ and
$-\div \, \sigma\,\mu=f$ on
$\ren\setminus \Sigma$. Then showing (\rf{conj1}) reduces to establish that
$$ \liminf_\e \int j^*(\sigma_\e) \, d\mu_\e\ \ge\ \int ({\bar j})^*(\sigma) \,
d\mu\ .\leqno(\frm{lowerbound})$$

Clearly we need more information on the oscillatory behaviour of
matrix fields
$\sigma_\e$ and we therefore consider a family of Young measures associated
with the triple $\{\sigma_\e,
\mu_\e,\mu\}$. It can be shown (see [9, Prop 4.3]) that there
exists a suitable family of
probability measures
$\nu_x$ on
$\renns$ (defined $\mu $ a.e.$x$ ) such that:
$$  \Psi(\sigma_\e)\, \mu_\e \ \rh \ \left(\int \Psi(\xi)\nu_x(d\xi)\right)
\, \mu \quad
\hbox{for all $\Psi \in C_0(\renns)$ \ . } \leqno(\frm{young})
$$
Testing (\rf{young}) in the particular case $\Psi(\xi)=\xi$ gives the equality
$\sigma = |\nu_x]$ between the weak limit $\sigma$ introduced above
and the barycenter of $\nu_x$.
Now our task is to provide arguments for establishing the following
inequality for every element
 $\nu$ of the Young family $\{\nu_x\}$ generated by $\{\sigma_\e\}$ (we
omit further the index
$x$) :
$$    \int j^* (\xi)\, d\nu \ge  ({\bar j})^*([\nu])  \ .\leqno(\frm{conj2})$$
Together with (\rf{young}), (\rf{conj2}) implies the lower bound inequality
(\rf{lowerbound}).
\medskip
A natural guess for proving (\rf{conj2}) is to conjecture that $\nu$ is
supported on the subset
 of  tensors with rank $<n$, yielding that $j^*(\xi)= ({\bar j})^*(\xi)$
holds $\nu$ a.e. This property
 suggests that if the sets  $A_\e$ oscillate at a small scale (so that for
example ${1_{A_\e}\over \e}$ converges weakly
to the Lebesgue measure), then  they need to concentrate at the same scale
on lower dimensional manifolds.
The reason for this is that $A_\e$ becomes thinner and thinner whereas, due
to the divergence condition, the unit
exterior normal $n_\e(x)$ satisfies
$\sigma_\e n_\e(x)=0$  on $\partial A_\e\setminus spt F$ .

Unfortunately this picture is not the good one as
we can see  on the following example in $\R^2$ suggested by P. Seppecher
[24]:
Take $\Om=Y:=(-1/2,1/2)^2$ and $r_\e$ chosen so that $\pi r_\e^2=\e$.
Denote by $\chi_\e$ the
$Y$-periodization of the characteristic function of the ball $\{
|y|<r_\e\}$. Then $\chi({x\over \e})$
determines a subset $A_\e$ of $\Om$ such that $|A_\e|\sim\e$ and
$\mu_\e={1_{A_\e}\over\e} dx \rh \mu$
being $\mu$  the Lebesgue measure on $\Om$. Now it is easy to construct  a
$Y$- periodic function
$\varphi_\e\in L^\infty( Y; \R^4_{sym})$ such that $\div \ \varphi_\e = 0$,
$\varphi_\e(y)= 0$ if $ |y|\ge 2\,r_\e$
and $\varphi_\e(y)= I_2 $ if $ |y|\le r_\e$ here $I_2$ denotes the identity
matrix).
Then the sequence $\sigma_\e= \varphi_\e({x\over\e})$ is divergence free
and satisfies (\rf{bound}).
However the Young family $\{\nu_x\}$  defined through (\rf{young}) turns
out to be independent on $x$ and
has a Dirac mass concentrated on the  matrix $I_2$.

 However we notice that in the previous example  the
weak limit of $\sigma_\e \mu_\e$ is zero so that (\rf{conj2}) still holds.
In fact the property
(\rf{conj2}) is far from requiring that the support of $\nu$ contains only
degenerate tensors. For example,
(\rf{conj2}) is satisfied if $\nu$ can be decomposed as a convex
combination of probability measures $\nu=
t \nu_0 + (1-t)
\nu_1$ where $\spt (\nu_0) \subset \{ \det \xi =0\}$ and $det ([\nu_1])=0$.
This can be easily checked
recalling that $j^* $ and ${\bar j}^*$ agree on determinant free tensors
and by making use of Jensen's
inequality.

Basically what we call ``{\it vanishing mass conjecture}'' consists in saying
that the property (\rf{conj2}) is
satisfied {\sl for every convex function $j$ on $\renns$ } and for all Young
measures generated by
sequences $\{\sigma_\e\}$ considered above. In particular we claim that the
validity of  (\rf{liminf}) or
(\rf{conj2}) is not related to the fact that $j$ is quadratic (although it
could be helpful to use tools like
compensated compactness or $H$-measures).
To conclude let us give   an intrinsic way to express (\rf{conj2})
independently of $j$ which has been suggested by
P.Seppecher [24]:
\medskip
\noindent {\bf Conjecture:} {\it The probability $\nu$ can be decomposed as
follows: there exists a
probability measure
$\nu_0$ on
$\{\det
\xi =0\}$ and for
$\nu_0$ almost all
$\xi$, there exists a probability measure $\lambda_\xi$ on $\renns$ such
that $[\lambda_\xi]=\xi$ and for all
$\Psi\in C_0(\renns)$:
$$   \int \Psi d\nu = \int \left(  \int \Psi(y) \lambda_\xi(dy)\right)
\nu_0(d\xi)\ .\leqno(\frm{conj3})$$ }
As before it is easy to check that (\rf{conj3}) implies (\rf{conj2}). Indeed by
using two times Jensen's inequality and the fact that $ j^*=({\bar j})^*$ \
 $\nu_0$ a.e., we obtain:
$$\eqalign{ \int j^* (\xi)\, d\nu &\ge \int \left(  \int j^*(y)
\lambda_\xi(dy)\right)\nu_0(d\xi) \ge
 \int j^*(\xi)\, \nu_0(d\xi) \cr
&\ge \int ({\bar j})^*(\xi)\, \nu_0(d\xi) \ge
({\bar j})^*([\nu])\
 .}$$


\bigskip
\noindent {\bf Acknowledgments.} I would like to thank warmly Giovanni Alberti and Ilaria Fragal\`a for
fruitful discussions on the topic and for sharing many attempts for
proving the {\it vanishing mass conjecture}. I thank also Pierre
Seppecher for suggesting me a counterexample and for proposing the
final statement (\rf{conj3}) as an alternative.


\bigskip
\noindent {\bf References.}

\item{[\bib{AK}]}
G. Allaire, R. V. Kohn: {\sl Optimal design for minimum
weight and compliance in plane stress using extremal microstructures.}
Europ. J. Mech. A/Solids, {\bf12} (6) (1993), 839--878.

\item{[\bib{BB}]}M. Bellieud, G. Bouchitt\'e: {\sl Homogenization of
elliptic problems in a fiber reinforced structure. Nonlocal effects.}
Ann. Scuola Norm. Sup. Pisa Cl. Sci., {\bf26} (1998), 407--436.

\item{[\bib{soon}]}
G. Bouchitt\'e: {\sl Optimal compliance of lower dimensional structures},
in preparation.

\item{[\bib{CV}]}
G. Bouchitt\'e, G. Buttazzo, P. Seppecher: {\sl Energies with
respect to a measure and applications to low dimensional structures.}
Calc. Var., {\bf 5} (1997), 37--54.

\item{[\bib{BoBu}]}
G. Bouchitt\'e, G. Buttazzo:
{\sl Characterization of optimal shapes and masses through
Monge-Kantorovich equation}.  J. Eur. Math. Soc. {\bf 3}
(2001), 139-168.

\item{[\bib{plap}]}
G. Bouchitt\'e, G. Buttazzo, G. De Pascale: {\sl A p-laplacian Approximation
for some Mass Optimization Problems}
preprint anam 2002-01, to appear in JOTA.

\item{[\bib{BoBuFr}]} G.~Bouchitt\'e, G.~Buttazzo, I.~Fragal\`a:
{\sl Convergence of Sobolev spaces on varying manifolds,} J.\
Geom.\ Anal.\ , Vol.~11 No.~3 (2001).

\item{[\bib{BoBuFr2}]} G.~Bouchitt\'e, G.~Buttazzo, I.~Fragal\`a:
{\sl Bounds for the effective coefficients of homogenized low-dimensional
structures} J.
Mat.\ Pures et Appl.\ , Vol.~81  (2002), 453 -- 469.
\item{[\bib{BoBuFra}]} G. Bouchitt\'e, G. Buttazzo and I. Fragal\`a,
{\sl Convergence of Sobolev spaces on varying manifolds,} to
Journal Geometrical Analysis, Vol.11, 3 ( 2001), 399--422.

\item{[\bib{BoBuSe}]}
G. Bouchitt\'e, G. Buttazzo, P. Seppecher:{\sl Shape optimization
solutions via Monge-Kantorovich equation.} C. R. Acad. Sci. Paris,
{\bf324-I} (1997), 1185--1191.
\item{[\bib{siam}]} G.\
Bouchitt\'e, I.\ Fragal\`a:{\sl Homogenization of thin structures by
two-scale method with respect to measures,} {\it SIAM J.\ Math.\
Anal.} {\bf 32} no.6 (2001), 1198--1226.
\item{[\bib{fgi}]} G.\ Bouchitt\'e, I.\ Fragal\`a:
{\sl Homogenization of elastic problems on thin structures: a
measure-fattening approach.} To appear in {\it J. Convex Anal.}

\item{[\bib{rev}]} G.\ Bouchitt\'e, I.\ Fragal\`a: {\sl Variational theory
of weak geometric structures},
Proceedings Como, to appear.

\item{[\bib{CaSe}]} M. Camar-Eddine, P.
       Seppecher. , {\sl Determination of the closure of the set of
elasticity functionals  },
preprint anam 2002/02

\item{[\bib{Cio}]}  D. Cioranescu and J. Saint Jean Paulin, Homogenization of
Reticulated Structures. Applied Mathematical Sciences {\bf 136},
Springer Verlag,  (1999).

\item{[\bib{FR}]} G. A. Francfort, F. Murat: {\sl Homogenization and
optimal bounds in linear elasticity.} Arch. Rational Mech. Anal., {\bf 94}
(1986), 307--334.

\item{[\bib{GS}]}
C. Goffmann, J. Serrin: {\sl Sublinear functions of measures and
variational integrals.}  Duke Math. J., {\bf31} (1964), 159--178.

\item{[\bib{GoS}]}
F. Golay, P. Seppecher: {\sl Locking materials and topology of optimal
shapes.}  Eur. J. Mech. A Solids {\bf 20} (4) (2001), 631--644.

\item{[\bib{KS}]}
R. V. Kohn, G. Strang: {\sl Optimal design and relaxation of
variational problems, I,II,III.} Comm. Pure Appl. Math., {\bf39}
(1986), p. 113-137, p. 139-182, 353--377.

\item{[\bib{M}]} U. Mosco: {\sl Composite media and asymptotic
Dirichlet forms.} J. Funct. Anal., {\bf123} (1994), 368--421.

\item{[\bib{Michell}]} A. Michell: {\sl The limits of economy of material
in frame-structures}, Phil. Mag.,8,
(1902)
589--597.

\item{[\bib{MT}]}
F. Murat, L. Tartar: {\sl Optimality conditions and homogenization.}
Proceedings of ``Nonlinear variational problems'', Isola d'Elba 1983,
Res. Notes in Math. {\bf127}, Pitman, London, (1985),  1--8.

\item{[\bib{MT1}]}
F. Murat, L. Tartar: {\sl Calcul des variations et
homog\'en\'eisation.} Proceedings of ``Les M\'ethodes de
l'homog\'en\'eisation: Th\'eorie et applications en physique'', Ecole
d'Et\'e d'Analyse Num\'erique C.E.A.-E.D.F.-INRIA,  Eyrolles, Paris,
(1985), 319--369.

\item{[\bib{Se}]} P. Seppecher: {\sl private communication}.

\item{[\bib{T}]}
L. Tartar: {\sl Estimations Fines des Coefficients Homog\'en\'eises.}
Ennio De Giorgi Colloqium, Edited by P.Kr\'ee, Res. Notes in Math.
{\bf125} Pitman, London (1985), 168--187.

\bigskip
Guy Bouchitt\'e \par
ANLA, UFR Sciences \par
Universit\'e de Toulon et du Var \par
BP 132, 83957 La Garde Cedex  \par
France \par
e-mail: bouchitte@univ-tl.fr
\medskip

\newwrite\outfile
\openout\outfile\jobname.out
\end